# Some Improved Estimators For Estimating Population Mean In Stratified Random Sampling

**Rajesh Singh, Viplav K. Singh, A. A. Adewara**

*Department of Statistics, Banaras Hindu University, Varanasi-221005, India

** Department of Statistics, University of Ilorin, Ilorin, Kwara State, Nigeria

**Abstract**

Some improved estimators are proposed for estimating the population mean in stratified sampling in the presence of auxiliary information. Mean square error (MSE) of the proposed estimators have been derived under large sample approximation. It has been shown that under optimum conditions proposed estimators are better than usual unbiased estimator and Hansen et al. (1946) estimator. Both theoretical and empirical findings are encouraging and support the soundness of the proposed procedure for mean estimation.

1. **Introduction**

Stratified sampling has often proved useful in planning surveys for improving the precision of other unstratified sampling strategies to estimate the finite population mean

$$\overline{Y} = \frac{1}{N} \sum_{h=1}^{L} \sum_{i=1}^{N_h} y_{hi.}$$

A ratio-product estimation of finite population mean $\overline{Y}$ can be made in two ways. One is to make a separate ratio-product estimate of the total of each stratum and add these totals. An alternative estimate is derived from a single combined ratio-product.

Consider a finite population of size N. Let y and x respectively, be the study and auxiliary variates on each unit $U_j$ ( j = 1,2,3,..., N ) of the population U. Let the population be divided into L strata with the $h^{th}$-stratum containing $N_h$ units, h=1,2,3,…,L so that $\sum_{h=1}^{L} N_h = N$.

Suppose that a simple random sample of size $n_h$ is drown without replacement from $h^{th}$ stratum such that $\sum_{h=1}^{L} n_h = n$.

We compute the sample mean of the variates in stratified sampling method as,

$$\bar{y}_{st} = \sum_{h=1}^{n_h} W_h \bar{y}_h \quad \text{and} \quad \bar{x}_{st} = \sum_{h=1}^{n_h} W_h \bar{x}_h$$

where,

$\bar{x}_h$ is the sample mean of auxiliary variates of $h^{th}$ stratum

$\bar{y}_h$ is the sample mean of study variates of $h^{th}$ stratum

$W_h = \dfrac{N_h}{N}$ is stratum weight.

The variance of usual unbiased estimator $\bar{y}_{st}$ is given as,

$$v(\bar{y}_{st}) = \sum_{h=1}^{N} w_h^2 \gamma_h S_{hy}^2$$

When the population mean $\bar{X}$ of the auxiliary variate x is known, Hansen, et al. (1946) suggested a "combined ratio estimator" as:

$$\bar{y}_{Rc} = \bar{y}_{st}\left(\dfrac{\bar{X}}{\bar{x}_{st}}\right) \tag{1.1}$$

The combined product estimator for $\bar{Y}$ is defined by,

$$\bar{y}_{Pc} = \bar{y}_{st}\left(\dfrac{\bar{x}_{st}}{\bar{X}}\right) \tag{1.2}$$

To the first degree of approximation, the variances of $\bar{y}_{Rc}$ and are $\bar{y}_{Pc}$ are respectively given by

$$V(\bar{y}_{Rc}) = \sum_{h=1}^{L} W_h^2 \gamma_h (S_{hy}^2 + R^2 S_{hx}^2 - 2RS_{hxy})  \qquad (1.3)$$

$$V(\bar{y}_{Pc}) = \sum_{h=1}^{L} W_h^2 \gamma_h (S_{hy}^2 + R^2 S_{hx}^2 + 2RS_{hxy})  \qquad (1.4)$$

where,

$$C_{hx}^2 = \frac{S_{hx}^2}{\overline{X}^2}, \quad C_{hy}^2 = \frac{S_{hy}^2}{\overline{Y}^2}, \quad R = \frac{\overline{Y}}{\overline{X}}, \quad \rho_{hxy} = S_{hxy} \frac{S_{hxy}}{S_{hx} S_{hy}}, \quad \gamma_h = \left(\frac{1}{n_h} - \frac{1}{N_h}\right)$$

$$S_{hy}^2 = \frac{1}{N_h - 1} \sum_{i=1}^{N_h} (y_{hi} - \overline{Y}_h)^2, \quad S_{hx}^2 = \frac{1}{N_h - 1} \sum_{i=1}^{N_h} (x_{hi} - \overline{X}_h)^2$$

$$S_{hxy} = \frac{1}{N_h - 1} \sum_{i=1}^{N_h} (x_{hi} - \overline{X}_h)(y_{hi} - \overline{Y}_h)$$

In this study, under stratified random sampling without replacement scheme, we suggest some improved estimators which are more efficient than estimator proposed by Hansen, et al. (1946) estimator.

**2. Proposed estimators**

Adapting Sahai and Ray (1980) estimator in stratified random sampling we propose an estimator $t_1$ as:

$$t_1 = \bar{y}_{st} \left[ 2 - \left\{ \frac{\bar{x}_{st}}{\overline{X}} \right\}^w \right]  \qquad (2.1)$$

We propose another estimator $t_2$ as:

$$t_2 = \bar{y}_{st} \left[ \frac{\bar{x}_{st} + a(\overline{X} - \bar{x}_{st})}{\bar{x}_{st} + b(\overline{X} - \bar{x}_{st})} \right]^p  \qquad (2.2)$$

To improve the efficiency of the estimators several authors have suggested combining ratio estimator with difference estimator in different ways. Some important references are Ray and Singh (1981), Singh et al. (2008), Gupta and Shabbir (2008), Grover and Kaur (2011) and Singh and Solanki (2012). Motivated by these authors we suggest some improved estimators combining ratio estimator with difference estimator as:

$$t_3 = \left( k_{31} \bar{y}_{st} + k_{32} (\overline{X} - \bar{x}_{st}) \right) \left[ 2 - \left\{ \frac{\bar{x}_{st}}{\overline{X}} \right\}^w \right]$$

(2.3)

$$t_4 = \left[k_{41}\bar{y}_{st} + k_{42}(\bar{X} - \bar{x}_{st})\right]\left(\frac{\bar{x}_{st} + a(\bar{X} - \bar{x}_{st})}{\bar{x}_{st} + b(\bar{X} - \bar{x}_{st})}\right)^P \quad (2.4)$$

$$t_5 = k_{51}\bar{y}_{st}\left(2 - \left(\frac{\bar{x}_{st}}{\bar{X}}\right)^w\right) + k_{52}(\bar{X} - \bar{x}_{st}) \quad (2.5)$$

$$t_6 = k_{61}\bar{y}_{st}\left[\frac{\bar{x}_{st} + a(\bar{X} - \bar{x}_{st})}{\bar{x}_{st} + b(\bar{X} - \bar{x}_{st})}\right]^P + k_{62}(\bar{X} - \bar{x}_{st}) \quad (2.6)$$

To obtain the biases and MSE's of the proposed estimators, we use the following notations in the rest of the article:

$$\bar{y}_{st} = \sum_{h=1}^{L} w_h \bar{y}_h = \bar{Y}(1 + e_0)$$

$$\bar{x}_{st} = \sum_{h=1}^{L} w_h \bar{x}_h = \bar{X}(1 + e_0)$$

Now expressing estimators in the terms of $e_i$'s (i=0,1), we have

$$t_1 = \bar{Y}\left[1 - we_1 - w(w-1)\frac{e_1^2}{2} + e_0 - we_0e_1\right] \quad (2.7)$$

$$t_2 = \bar{Y}\left[1 + \{e_0 + e_1 D\} + e_1^2 C + e_0 e_1\right] \quad (2.8)$$

$$t_3 = k_{31}\bar{Y}e_0 + k_{31} - wk_{31}e_0e_1\bar{Y} - wk_{31}e_1\bar{Y} - \frac{w(w-1)}{2}k_{31}e_1^2\bar{Y} - k_{32}e_1\bar{X} + wk_{32}e_1^2\bar{X} \quad (2.9)$$

$$t_4 = \bar{Y}k_{41}[1 + e_0 - pe_1(1-b) + pe_1(1-a)] + k_{41}\bar{Y}e_1^2 \frac{[p^2(a-b)^2 + p(b^2 - a^2 + 2ab)]}{2}$$
$$+ k_{41}\bar{Y}e_0e_1p(b-a) - k_{42}e_1\bar{X} + k_{42}pX e_1^2(1-b) - k_{42}p\bar{X}(1-a)e_1^2 \quad (2.10)$$

$$t_5 = k_{51}\bar{Y}[1 + e_0 - we_0e_1 - we_1 - w(w-1)\frac{e_1^2}{2}] - k_{52}\bar{X}e_1 \quad (2.11)$$

$$t_6 = k_{61}Y[1 + (E_0 + p(b-a)e_1)] + e_{61}^2 C - Ae_0e_1 - k_{62}\bar{X}e_1 \quad (2.12)$$

Taking expectations and than substracting $\bar{Y}$, we get the biases of the above estimators,

respectively as:

$$B(t_1) = \overline{Y}\left[\frac{w(1-w)}{2}\frac{v(\overline{x}_{st})}{\overline{X}^2} - \frac{w\,cov(\overline{y}_{st},\overline{x}_{st})}{\overline{Y}\overline{X}}\right] \tag{2.13}$$

$$B(t_2) = \overline{Y}\left[C\frac{v(\overline{x})}{\overline{X}} - A\frac{cov(\overline{x}_{st},\overline{y}_{st})}{\overline{X}\overline{Y}}\right] \tag{2.14}$$

$$B(t_3) = \overline{Y}(k_{31}-1) - k_{31}\overline{Y}\left\{\frac{w(w-1)v(\overline{x}_{st})}{2\overline{X}^2} + \frac{w\,cov(\overline{x}_{st}\overline{y}_{st})}{\overline{X}\overline{Y}}\right\} + k_{32}\overline{X}w\,\frac{v(\overline{x}_{st})}{\overline{X}^2} \tag{2.15}$$

$$B(t_4) = \overline{Y}(k_{41}-1) + k_{41}\overline{Y}\left\{\frac{v(\overline{x})}{\overline{X}^2}C + D'\frac{cov(\overline{x}_{st},\overline{y}_{st})}{\overline{X}\overline{Y}}\right\} - k_{42}D'\frac{v(\overline{x}_{st})}{\overline{X}} \tag{2.16}$$

$$B(t_5) = \overline{Y}(k_{51}-1) + k_{51}\overline{Y}\left\{\frac{w\,cov(\overline{x}_{st}\overline{y}_{st})}{\overline{X}\overline{Y}} - w(w-1)\frac{v(\overline{x}_{st})}{2\overline{X}^2}\right\} \tag{2.17}$$

$$B(t_6) = \overline{Y}(k_{61}-1) + k_{61}\overline{Y}\left\{\frac{v(\overline{x}_{st})}{\overline{X}^2}C - A\frac{cov(\overline{x}_{st}\overline{y}_{st})}{\overline{X}\overline{Y}}\right\} \tag{2.18}$$

where,

$$A = p(1-a)$$
$$B = p(1-b)$$
$$C = \frac{p^2(a-b)^2 + p[b^2 - a^2 + 2(a-b)]}{2}$$
$$D' = p(a-b).$$

The MSE expressions of the above estimator's are respectively given by

$$MSE(t_1) = v(\overline{y}_{st}) + w^2 R^2 v(\overline{x}_{st}) - 2wR\,cov(\overline{y}_{st},\overline{x}_{st}) \tag{2.19}$$

$$MSE(t_2) = v(\overline{y}_{st}) + D'^2 R^2 v(\overline{x}_{st}) + 2D'R\,cov(\overline{y}_{st},\overline{x}_{st}) \tag{2.20}$$

$$MSE(t_3) = \overline{Y}^2(K_{31}-1)^2 + K_{31}^2\left\{v(\overline{y}_{st}) + w\,R^2 v(\vec{x}_{st}) - 4wR\,cov(\overline{y}_{st},\vec{x}_{st})\right\}$$
$$+ K_{32}^2 v(\vec{x}_{st}) + 2K_{31}\left\{wR\,cov(\overline{y}_{st},\vec{x}_{st}) + \frac{w(w-1)R^2 v(\overline{x}_{st})}{2}\right\} - 2K_{32}wRv(\overline{x}_{st})$$
$$+ 2K_{31}K_{32}\left\{2wRv(\vec{x}_{st}) - cov(\overline{y}_{st},\vec{x}_{st})\right\} \tag{2.21}$$

$$\text{MSE}(t_4) = \overline{Y}^2(K_{41}-1)^2 + K_{41}^{\ 2}A_4 + K_{42}^{\ 2}B_4 - 2K_{41}C_4 + 2K_{42}D_4 - 2K_{41}K_{42}E_4 \qquad (2.22)$$

$$\text{MSE}(t_5) = \overline{Y}^2(K_{51}-1)^2 + K_{51}^{\ 2}\left\{v(\overline{y}_{st}) + w^2R^2v(\vec{x}_{st}) - 4wR\,\text{cov}(\overline{y}_{st},\vec{x}_{st}) - w(w-1)R^2v(\vec{x}_{st})\right\}$$
$$+ K_{52}^{\ 2}v(\vec{x}_{st}) + 2K_{51}\left\{wR\,\text{cov}(\overline{y}_{st},\vec{x}_{st}) + \frac{w(w-1)R^2v(\vec{x}_{st})}{2}\right\}$$
$$- 2K_{51}K_{52}\left\{\text{cov}(\overline{y}_{st},\vec{x}_{st}) - wRv(\vec{x}_{st})\right\} \qquad (2.23)$$

$$\text{MSE}(t_6) = \overline{Y}^2(K_{61}-1)^2 + K_{61}^{\ 2}\left\{v(\overline{y}_{st}) + D'^2R^2v(\vec{x}_{st}) + 4D'R\,\text{cov}(\overline{y}_{st},\vec{x}_{st}) + 2R^2v(\vec{x}_{st})C\right\}$$
$$+ K_{62}^{\ 2}v(\vec{x}_{st}) - 2K_{61}\left\{R^2v(\vec{x}_{st})C + RD'\text{cov}(\overline{y}_{st},\vec{x}_{st})\right\}$$
$$- 2K_{61}K_{62}\left\{\text{cov}(\overline{y}_{st},\vec{x}_{st}) + D'Rv(\vec{x}_{st})\right\} \qquad (2.24)$$

Where,

$$A_4 = \overline{Y}^2\left(\frac{v(\overline{y}_{st})}{\overline{Y}^2} + D'^2\frac{v(\overline{x}_{st})}{\overline{X}^2} + 4D'\frac{\text{cov}(\overline{y}_{st},\overline{x}_{st})}{\overline{X}\overline{Y}} + 2C\frac{v(\overline{y}_{st})}{\overline{X}^2}\right)$$

$$B_4 = \overline{X}^2\frac{v(\overline{x}_{st})}{\overline{X}^2}$$

$$C_4 = \overline{Y}^2\left(C\frac{v(\overline{x}_{st})}{\overline{X}^2} + D'\frac{\text{cov}(\overline{y}_{st},\overline{x}_{st})}{\overline{X}\overline{Y}}\right)$$

$$D_4 = \overline{X}\frac{v(\overline{x}_{st})}{\overline{X}^2}D'$$

$$E_4 = \overline{X}\overline{Y}\left(D'\frac{v(\overline{x}_{st})}{\overline{X}^2} + \frac{\text{cov}(\overline{y}_{st},\overline{x}_{st})}{\overline{X}\overline{Y}}\right)$$

Partially differentiating equation (2.21) with respect to $K_{31}$ and $K_{32}$, we get the optimum values as:

$$K_{31}(\text{opt}) = \frac{D_3E_3 - B_3(\overline{Y}^2 - C_3)}{E_3^{\ 2} - B_3(\overline{Y}^2 + A_3)}, \qquad K_{32}(\text{opt}) = \frac{E_3(\overline{Y}^2 - C_3) - D_3(\overline{Y}^2 + A_3)}{E_3^{\ 2} - B_3(\overline{Y}^2 + A_3)}$$

Where,

$$A_3 = \left\{v(\overline{y}_{st}) + w\,R^2v(\vec{x}_{st}) - 4wR\,\text{cov}(\overline{y}_{st},\vec{x}_{st})\right\}$$

$$B_3 = v(\vec{x}_{st})$$

$$C_3 = \left\{wR\,\text{cov}(\overline{y}_{st},\vec{x}_{st}) + \frac{w(w-1)R^2v(\vec{x}_{st})}{2}\right\}$$

$$D_3 = wRv(\vec{x}_{st})$$

$$E_3 = \{2wRv(\vec{x}_{st}) - cov(\bar{y}_{st}, \vec{x}_{st})\}$$

Similarly, partially differentiating equation (2.22) with respect to $K_{41}$ and $K_{42}$, we get the optimum values as:

$$k_{41}(opt) = \frac{B_4(\bar{Y}^2 + C_4) - D_4 E_4}{B_4(\bar{Y}^2 + A_4) - E_4^2} \ , \ k_{42}(opt) = \frac{E_4(\bar{Y}^2 + C_4) - D_4(\bar{Y}^2 + A_4)}{B_4(\bar{Y}^2 + A_4) - E_4^2} \quad (2.26)$$

Where,

$$A_4 = \bar{Y}^2 \left( \frac{v(\bar{y}_{st})}{\bar{Y}^2} + D'^2 \frac{v(\bar{x}_{st})}{\bar{X}^2} + 4D' \frac{cov(\bar{y}_{st}, \bar{x}_{st})}{\bar{X}\bar{Y}} + 2C \frac{v(\bar{y}_{st})}{\bar{X}^2} \right)$$

$$B_4 = \bar{X}^2 \frac{v(\bar{x}_{st})}{\bar{X}^2}$$

$$C_4 = \bar{Y}^2 \left( C \frac{v(\bar{x}_{st})}{\bar{X}^2} + D' \frac{cov(\bar{y}_{st}, \bar{x}_{st})}{\bar{X}\bar{Y}} \right)$$

$$D_4 = \bar{X} \frac{v(\bar{x}_{st})}{\bar{X}^2} D'$$

$$E_4 = \bar{X}\bar{Y} \left( D' \frac{v(\bar{x}_{st})}{\bar{X}^2} + \frac{cov(\bar{y}_{st}, \bar{x}_{st})}{\bar{X}\bar{Y}} \right)$$

Now, partially differentiating equation (2.23) with respect to $K_{51}$ and $K_{52}$, we get the optimum values as:

$$k_{51}(opt) = \frac{B_5(\bar{Y}^2 - C_5)}{B_5(\bar{Y}^2 + A_5)D_5^2} \quad K_{52}(opt) = \frac{D_5(\bar{Y}^2 - C_5)}{B_5(\bar{Y}^2 + A_5)D_5^2} \quad (2.27)$$

Where,

$$A_5 = \{v(\bar{y}_{st}) + w^2 R^2 v(\vec{x}_{st}) - 4wR \, cov(\bar{y}_{st}, \vec{x}_{st}) - w(w-1)R^2 v(\vec{x}_{st})\}$$

$$B_5 = v(\vec{x}_{st})$$

$$C_5 = \left\{ wR \, cov(\bar{y}_{st}, \vec{x}_{st}) + \frac{w(w-1)R^2 v(\vec{x}_{st})}{2} \right\}$$

$$D_5 = \{cov(\bar{y}_{st}, \vec{x}_{st}) - wRv(\vec{x}_{st})\}$$

Finally, partially differentiating equation (2.24) with respect to $K_{51}$ and $K_{52}$, we get the optimum values as:

$$k_{61}(opt) = \frac{(\bar{Y}^2 + C_6)B_6}{B_6(\bar{Y}^2 + A_6) - D_6^2} \quad , \quad k = \frac{(\bar{Y}^2 + C_6)D_6}{B_6(\bar{Y}^2 + A_6) - D_6^2} \tag{2.28}$$

Where,

$$A_6 = \{v(\bar{y}_{st}) + D'^2 R^2 v(\vec{x}_{st}) + 4D'R\, cov(\bar{y}_{st}, \vec{x}_{st}) + 2R^2 v(\vec{x}_{st})C\}$$

$$B_6 = v(\vec{x}_{st})$$

$$C_6 = \{R^2 v(\vec{x}_{st})C + RD'cov(\bar{y}_{st}, \vec{x}_{st})\}$$

$$D_6 = \{cov(\bar{y}_{st}, \vec{x}_{st}) + DRv(\vec{x}_{st})\}$$

## 3. Empirical Study

To see the performance of various estimators of population mean $\bar{Y}$, with respect to Usual unbiased estimator $\bar{y}_{st}$, we have considered two data sets. Summaries of the Data are given below:

**Data set 1:** Source Singh and Mangat:

$y_h$ : Juice quantity, $\quad x_h$ : Weight of cane.

| Total | Stratum | 1 | 2 | 3 |
|---|---|---|---|---|
| N=25 | $N_h$ | 6 | 12 | 7 |
| n=10 | $n_h$ | 3 | 4 | 3 |
| $\overline{X}$ =326 | $\overline{X}_h$ | 366.666 | 310.883 | 317.143 |
| $\overline{Y}$ =102.6 | $\overline{Y}_h$ | 135 | 99.166 | 80.714 |
| $S_x^2$ =2700 | $S_{xh}^2$ | 2706.666 | 1881.06 | 2890.476 |
| $S_y^2$ =558.583 | $S_{yh}^2$ | 80 | 226.515 | 120.238 |
| $\rho$ =.7314955 | $\rho_h$ | 0.9455626 | 0.948196 | 0.7523324 |
| R=0.314723 | $\gamma_h$ | 0.1666667 | 0.1666667 | 0.1904762 |
| $\rho_c$ =0.8676778 | $W_h^2$ | 0.0576 | 0.2304 | 0.0784 |

**DATA SET 2:** Source Singh and Chaudhary (1986, pg.162)

The data were collected in a pilot survey for estimating the extent of cultivation and production of fresh fruits in three districts of U.P .

$x_h$ : area under orchards in hect.

$y_h$ : total no of trees

| Stratum | 1 | 2 | 3 |
|---|---|---|---|
| $N_h$ | 985 | 2196 | 1020 |
| $n_h$ | 6 | 8 | 11 |
| $\overline{X}_h$ | 11253 | 25115 | 18870 |
| $S_{xh}^2$ | 15.97 | 132.66 | 38.44 |
| $S_{yh}^2$ | 74775.47 | 259113.7 | 65885.6 |
| $S_{xyh}$ | 1007.75 | 5709.16 | 1404.71 |
| $\gamma_h$ | 0.16598 | 0.12454 | 0.08902 |

and
$R = 49.03, \alpha_{opt} = 0.9422$

# Table- 3.1: MSE's and PRE's of estimators

|   | ESTIMATORS | Data-1 | | Data-2 | |
|---|---|---|---|---|---|
|   |   | MSE | PRE | MSE | PRE |
| 1 | $t_1$ | 701.546 | 1403.318 | 2.782946 | 404.6695 |
| 2 | $t_2$ | 701.54 | 1403.318 | 2.782946 | 404.6695 |
| 3 | $t_3$ | 629.0631 | 1565.013 | 2.77094 | 404.9483 |
| 4 | $t_4$ | 874.5025 | 1125.774 | 3.051538 | 369.0511 |
| 5 | $t_5$ | 601.846 | 1635.864 | 2.77668 | 405.5826 |
| 6 | $t_6$ | 524.6948 | 1876.314 | 2.77092 | 406.4257 |
| 7 | $\bar{y}_{RC}$ | 857.37974 | 1148.2567 | 3.47243 | 324.3185 |
| 8 | $\bar{y}_{PC}$ | 21953.129 | 44.84 | 47.0589 | 23.93111 |
| 9 | $\bar{y}_{st}$ | 9844.9203 | 100.00 | 11.26173 | 100.00 |

**Conclusion**

From Table 3.1, we see that all the proposed estimators perform better than usual mean estimator and combined ratio estimator. For data set 1 estimator $t_6$ is best followed by the estimators $t_5$ and $t_3$. For data set 2 $t_6$ is the best estimator followed by the estimator $t_5$.